\documentclass[10pt,twoside]{amsart}
\usepackage{lmodern}
\usepackage[T1]{fontenc}
\usepackage{amsmath,amsfonts,amsthm,amssymb,amscd}
\usepackage{mathrsfs}
\usepackage{euscript}
\usepackage{latexsym}
\usepackage[pdftex]{graphicx}
\usepackage{enumitem}
\usepackage{bbm}
\usepackage{accents}

\usepackage{geometry}
\newtheorem{corollary}{Corollary}[section]
\newtheorem{lemma}{Lemma}[section]
\newtheorem{theorem}{Theorem}[section]

\usepackage{hyperref}  %added
\hypersetup{colorlinks=true, linktoc=all,
    linkcolor=blue,
		citecolor=blue,
		urlcolor=cyan,
		bookmarksopen=true
		}
\usepackage{hyperref}  %added
\newcommand{\lnc}{\mathscr{L}}

\title{On the mean values of the Chebyshev functions and their applications}
\author{Zarullo Rakhmonov, Opokkhon Nozirov}
\address{A.Dzhuraev Institute of Mathematics,  National Academy of Sciences of Tajikistan}
\email{zarullo.rakhmonov@gmail.com, zarullo-r@rambler.ru, nozirov92@inbox.ru}
\date{}

\begin{document}

\begin{abstract}
Assuming the validity of the extended Riemann hypothesis for the average values of Chebyshev functions over all characters modulo $q$, the following estimate holds
$$
t(x;q)=\sum_{\chi\bmod q}\max_{y\leq x}|\psi(y,\chi)|\ll x+x^\frac12q{\mathscr{L}}^2,\qquad \mathscr{L}=\ln xq.
$$
When solving a number of problems in prime number theory, it is sufficient that $t(x;q)$ admits an estimate close to this one. The best known estimates for $t(x;q)$ previously belonged  to G.~Montgomery, R.~Vaughn, and Z.~Kh.~Rakhmonov. In this paper we obtain a new estimate of the form 
$$
t(x;q)=\sum_{\chi\bmod q}\max_{y\leq x}|\psi(y,\chi)|\ll x{\mathscr{L}}^{28}+x^\frac45q^\frac12{\mathscr{L}}^{31}+x^\frac12q{\mathscr{L}}^{32},
$$
using which for a linear exponential sum with primes we prove a stronger estimate
$$
S(\alpha,x)\ll xq^{-\frac12}\lnc^{33}+x^{\frac{4}{5}}\lnc^{32}+x^\frac{1}{2}q^\frac12\lnc^{33},
$$
when $\left|\alpha-\frac aq\right|<\frac1{q^2}$,  $(a,q)=1$. We also study the distribution of Hardy-Littlewood numbers of the form $ p + n ^ 2 $ in short arithmetic progressions in the case when the difference of the progression is a power of the prime number.
Bibliography: 30 references.
\end{abstract}

\maketitle

\section{Introduction}
For a Dirichlet character $\chi$ modulo $q$, the Chebyshev function is defined by the equality
$$
\psi(y,\chi)=\sum_{n\le y}\Lambda(n)\chi(n),
$$
where $\Lambda(n)$ is the von Mangoldt function. Assuming the validity of the extended Riemann hypothesis for the mean values of Chebyshev functions over all characters of modulus $q$, the following estimate holds:
\begin{align}
t(x;q)=\sum_{\chi\hspace{-9pt} \ \mod\hspace{-2pt}q}\max_{y\leq x}|\psi(y,\chi)|\ll x+x^{1/2}q\lnc^2,\qquad \lnc=\ln xq.\label{formula-otsenka t(x;q) v pred RGR}
\end{align}
For solving a number of problems in prime number theory, it suffices that $t(x;q)$ admits an estimate close to (\ref{formula-otsenka t(x;q) v pred RGR}).

The study of mean values of Chebyshev functions was first undertaken by Yu.~V.~Linnik \cite{Linnik-1980-Nauka,Linnik-1946-MatSbor,Linnik-1945-DANSSSr,Linnik-1946-IzvANSSSr} in order to derive a nontrivial estimate for the linear exponential sum with prime numbers $S(\alpha,x)$.

A.~A.~Karatsuba \cite{Karatsuba-DANSSSR-1970-192-4} developed a method for solving ternary multiplicative problems, which he used to estimate the simplest case of $t(x;q)$. As a consequence of this estimate, the distribution of numbers of the form $p(p_1+a)$ in short arithmetic progressions was obtained.

Using Linnik`s large sieve method, G.~Montgomery \cite{Montgomeri-1974} proved density theorems for the zeros of Dirichlet $L$-functions, which allowed him to show that
\begin{align}\label{formula otcenka Montgomeri dlya t(x;q)}
&t(x;q)\ll(x+x^{\frac57}q^{\frac57}+x^{\frac12}q)\lnc^{17}.
\end{align}

This result was refined by R.~Vaughan \cite{Vaughan-1975}, who, using the representation
\begin{align*}
\frac{L'}{L}=\left(\frac{L'}{L}+F\right)(1-FG)+(L'+LF)G-F,
\end{align*}
where $F$ and $G$ are partial sums of the Dirichlet series $\frac{L'}{L}$ and $\frac{1}{L}$, respectively, proved that
\begin{align}\label{formula otcenka Vaughan dlya t(x;q)}
&t(x;q)\ll x\lnc^3+x^{\frac34}q^{\frac58}\lnc^{\frac{23}{8}}+x^{\frac12}q\lnc^{\frac72}.
\end{align}
In 1989, Z.~Kh.~Rakhmonov \cite{RakhmonovZKh-Izv-ANSSSR-1989} showed that
$$
t(x;q)\ll(x+x^{\frac56}q^{\frac12}+x^{\frac12}q)x^\delta.
$$
This estimate is stronger than (\ref{formula otcenka Montgomeri dlya t(x;q)}) but weaker than (\ref{formula otcenka Vaughan dlya t(x;q)}). However, unlike these estimates, its proof is elementary and is based on A.~A.~Karatsuba's method for solving ternary multiplicative problems \cite{Karatsuba-DANSSSR-1970-192-4}.

From the estimates (\ref{formula-otsenka t(x;q) v pred RGR}), (\ref{formula otcenka Montgomeri dlya t(x;q)}), and (\ref{formula otcenka Vaughan dlya t(x;q)}) for $t(x;q)$, it follows that among the three terms present in these bounds, the first and the last are equal up to a finite power of the logarithm, and they apparently cannot be improved in terms of the exponents of $x$ and $q$.

Further improvement of the second term was achieved by Z.~Kh.~Rakhmonov \cite{RakhmonovZKh-1993-IzvRAN,RakhmonovZKh-DAN-1993}, who proved that
\begin{align}\label{formula Rakhmonova-1993 otcenka dlya t(x;q)}
t(x;q)\ll \left(x+x^{\frac{4}{5}}q^{\frac12}+x^\frac{1}{2}q\right)\lnc^{34}.
\end{align}

The following theorem refines this estimate.
{\theorem\label{Teorema o srednem znach funk Chebisheva} For $x\ge 2$ and $q\ge 1$, the following estimate holds:
\begin{align*}
&t(x;q)\ll x\lnc^{28}+x^{\frac{4}{5}}q^{\frac12}\lnc^{31}+x^\frac{1}{2}q\lnc^{32}.
\end{align*}}

In 1937, I.~M.~Vinogradov \cite{VinigradovIM-Izb-Trud} discovered that sums over prime numbers can be expressed solely through addition and subtraction of a relatively small number of other sums, for which good estimates can be obtained using the method of bounding double sums, independent of the theory of Dirichlet $L$-series. In particular, such a sum is the linear exponential sum with prime numbers of the form
$$
S(\alpha,x)=\sum_{\substack{n\le x}}\Lambda(n)e(\alpha n),
$$
where $\alpha$ is a real number, and under the condition
$$
\left|\alpha-\frac aq\right|\leq \frac1{q^2},\qquad q\le x,\qquad (a,q)=1,
$$
he established the estimate:
\begin{align}\label{formula otsenka IMV}
&S(\alpha,x)\ll (xq^{-\frac12}+x^{\frac45}+x^{\frac12}q^{\frac12})x^{\varepsilon},
\end{align}
whose proof is based on an elementary method.

The sum $S(\alpha,x)$ was first estimated using an analytic method by Yu.~V.~Linnik \cite{Linnik-1980-Nauka,Linnik-1946-IzvANSSSr} (see also \cite{Chudakov-1947-Math,Chudakov-1947}). Using ideas of Hardy and Littlewood \cite{Hardy+Littlwood-1923}, previously applied in the Goldbach problem, along with the density theorem for the zeros of Dirichlet $L$-series, he provided a new nontrivial bound for the linear exponential sum with prime numbers in the following form: \emph{let $\alpha$ be a real number, $N\geq N_0>0$, $\alpha=\frac aq+\lambda$, where $(a,q)=1$, $1<q\le \tau=(\ln x)^{1000}$, $\tau^{1000}x^{-1}\le|\lambda|\le(q\tau)^{-1}$, then the estimate holds}
\begin{align*}
&\left|\sum_{\substack{n=2}}^{\infty}\Lambda(n)\exp\left(-\frac nN\right)e(\alpha n)\right|<N(\ln N)^{-1000}.
\end{align*}

G.~Montgomery \cite{Montgomeri-1974}, using his estimate for the mean values of Chebyshev functions (\ref{formula otcenka Montgomeri dlya t(x;q)}), proved that
\begin{align}\label{formula otcenka Montgomeri S(aq,x)}
&S\left(\frac aq,x\right)\ll\left(xq^{-\frac12}+x^{\frac57}q^{\frac{3}{14}}+x^{\frac12}q^{\frac12}\right)\lnc^{17}.
\end{align}
He also showed that if $\eta\le x^{\frac14}$, $\eta\le q\le x\eta^{-1}$, and $|\alpha-a/q|\le2\eta(qx)^{-1}$, where $(a,q)=1$, then
\begin{align}\label{formula otcenka Montgomeri S(alpha,x)-eta}
&S(\alpha,x)\ll x\eta^{-\frac12}\lnc^{17}.
\end{align}

R.~Vaughan \cite{Vaughan-1975}, applying his estimate for the mean values of Chebyshev functions (\ref{formula otcenka Vaughan dlya t(x;q)}), refined Montgomery's result. He proved that if $|\alpha-a/q|\le q^{-2}$, $(a,q)=1$, then the following estimate holds:
\begin{align}
&S(\alpha,x)\ll(xq^{-\frac12}+x^\frac78q^{-\frac18}+x^\frac34q^\frac18+x^\frac12q^\frac12)\lnc^4,\label{formula otsenka Vaughan S(alpha,x)}
\end{align}
and if $\eta\le x^{\frac13}$, $\eta\le q\le x\eta^{-1}$, and $|\alpha-a/q|\le2\eta(qx)^{-1}$, where $(a,q)=1$, then
\begin{align}\label{formula otsenka Vaughan S(alpha,x)-eta}
&S(\alpha,x)\ll x\eta^{-\frac12}\lnc^{4}.
\end{align}

It is worth noting that the estimates (\ref{formula otcenka Montgomeri S(aq,x)}), (\ref{formula otcenka Montgomeri S(alpha,x)-eta}), (\ref{formula otsenka Vaughan S(alpha,x)}), and (\ref{formula otsenka Vaughan S(alpha,x)-eta}), obtained by analytic methods, are weaker than the estimate (\ref{formula otsenka IMV}), which was obtained by I.~M.~Vinogradov using an elementary method.

Z.~Kh.~Rakhmonov \cite{RakhmonovZKh-1993-IzvRAN,RakhmonovZKh-DAN-1993}, utilizing his estimate for the mean values of Chebyshev functions (\ref{formula Rakhmonova-1993 otcenka dlya t(x;q)}), derived an estimate in which the factor $x^\varepsilon$ in (\ref{formula otsenka IMV}) is replaced by a finite power of the logarithm. That is, if $|\alpha-a/q|\le q^{-2}$, $(a,q)=1$, then
\begin{equation}\label{formula otsenka Rakhmonov S(alpha,x)}
S(\alpha,x)\ll(xq^{-\frac12}+x^{\frac45}+x^{\frac12}q^{\frac12})\lnc^{35}.
\end{equation}
Moreover, if $1\le\eta\le x^\frac25$, $\eta\le q<x\eta^{-1}$, and $|\alpha-a/q|\le2\eta(qx)^{-1}$, where $(a,q)=1$, then the following estimate holds:
\begin{align}\label{formula otsenka Rakhmonov S(alpha,x)-eta}
&S(\alpha,x)\ll x\eta^{-\frac12}\lnc^{35}.
\end{align}

Utilizing Theorem \ref{Teorema o srednem znach funk Chebisheva}, we establish an estimate for the sum $S(\alpha,x)$, refining the logarithmic exponents in the terms of (\ref{formula otsenka Rakhmonov S(alpha,x)}). Specifically, the estimate (\ref{formula otsenka Rakhmonov S(alpha,x)}) is first refined for the case when $\alpha$ is a rational number (Theorem \ref{Teorema otsenka S(a/q,x)}) and subsequently for an arbitrary real number $\alpha$ (Corollary \ref{Sledstvie1 teoremi otsenka S(alpha,x)}).

{\theorem\label{Teorema otsenka S(a/q,x)} Let $(a,q)=1$. Then the following estimate holds:
\begin{align*}
S\left(\frac aq,x\right)\ll xq^{-\frac12}\lnc^{29}+x^{\frac{4}{5}}\lnc^{32}+x^\frac{1}{2}q^\frac12\lnc^{33}.
\end{align*}}

{\corollary\label{Sledstvie1 teoremi otsenka S(alpha,x)} Let $\left|\alpha-\dfrac aq\right|<\dfrac{1}{q^2}$, \ $(a,q)=1$, then the following estimate holds:
$$
S(\alpha,x)\ll xq^{-\frac12}\lnc^{33}+x^{\frac{4}{5}}\lnc^{32}+x^\frac{1}{2}q^\frac12\lnc^{33}.
$$}

The next corollary provides a refinement of the estimate (\ref{formula otsenka Rakhmonov S(alpha,x)-eta}).
{\corollary\label{Sledstvie2 teoremi otsenka S(alpha,x)} Let $\eta\le x^\frac25$, $\eta\le q\le x\eta^{-1}$, $|\alpha-aq^{-1}|\le 2\eta(qx)^{-1}$, $(a,q)=1$, then the following estimate holds:
$$
S(\alpha,x)\ll x\eta^{-\frac12}\lnc^{33}.
$$}

Hardy and Littlewood \cite{Hardy+Wright} formulated the hypothesis that all sufficiently large natural numbers $n$ can be represented as the sum of a prime and a power of a natural number:
$$
n=p+m^{k},\qquad k\geq 2.
$$
Numbers satisfying this representation are referred to as Hardy-Littlewood numbers. G.~Babaev \cite{Babaev-UMN-1958} disproved this hypothesis by demonstrating the existence of infinitely many natural numbers that are not Hardy-Littlewood numbers. Consequently, it follows that there exist values of $l$, where $1\leq l\leq q$, such that the inequality
$$
H_{k}(q,l)>q, \qquad k\geq 2,
$$
holds, where $H_{k}(q,l)$ is the smallest Hardy-Littlewood number of the form $p+m^{k}$ within the arithmetic progression $qt+l$, $t=0,1,2,\ldots$, with $q$ being an integer. Therefore, it is natural to consider the following two problems:
\begin{enumerate}
  \item To obtain an upper bound for $H_{k}(q,l)$ with the best possible precision.
  \item To derive an asymptotic distribution law for Hardy-Littlewood numbers in very short arithmetic progressions.
\end{enumerate}

For the case where $q$ is a prime number and $k\geq 2$, these problems have been investigated in \cite{RakhmonovZKh-Izv-ANSSSR-1989,RakhmonovZKh-1993-IzvRAN,RakhmonovZKh-DAN-1993,RakhmonovZKh-1994-TrMIRAN}, yielding an asymptotic formula for the number of solutions to the congruence:
\begin{align*}
&p+m^{k}\equiv l\hspace{-8pt}\pmod q,\qquad   p\leq x, \qquad  m\leq \sqrt[k]{x}, \\
&q\ll\min\left(x^{\frac2k}\lnc_x^{-8},x^{\frac{k+5}{5k}}\lnc_x^{-35},x^{\frac{k+2}{3k}}\lnc_x^{-\frac{70}{3}},\right),\qquad \lnc_x=\ln x,
\end{align*}
which, in particular, implies that
$$
H_{2}(q,l)\ll q^\frac{3}{2} \ln^{35}q.
$$

The proof of this result is based on A.~A.~Karatsuba's method for solving multiplicative ternary problems \cite{Karatsuba-DANSSSR-1970-192-4} and on A.~Weil's theorem \cite{Weil-A-1948} for estimating complete mixed sums of the form
$$
S(\chi,g,f,p^\beta)=\sum_{m=1}^{p^\beta}\chi(g(m))e\left(\frac{f(m)}{p^\beta}\right),
$$
for the case $\beta=1$, where $\chi$ is a Dirichlet character modulo $p^\beta$, and $g(m)$ and $f(m)$ are rational functions defined modulo $p^\beta$, with $g(m)$ being nonzero modulo $p$.

The following theorem generalizes and refines this result for the case $k=2$ and when $q$ - the difference in the progression, is a power of a prime number.

{\theorem\label{Teorema o raspr chisel Hardy-Litlwoda} Let $x\ge x_0$, $p$ be an odd prime, $(l,p)=1$, and $-l$ be a quadratic non-residue modulo $p$. Define $\rho(p,l)$ as the number of solutions to the congruence $n^2\equiv l\hspace{-5pt}\pmod p$,
$$
\EuScript{H}_2(x;p^\alpha,l)=\sum_{\substack{n\le x,\ m^2\leq x  \\ n+m^2\equiv l\hspace{-7pt}\pmod{p^\alpha}}}\Lambda(n).
$$
Then for any fixed $A\ge 58$, the following asymptotic formula holds:
$$
\EuScript{H}_2(x;p^\alpha,l)=\frac{x^\frac32}{\varphi(p^\alpha)}\left(1-\frac{\rho(p,l)}{p}+O\left(\lnc_x^{-0.5A+28}+ \frac{p^{0.5\alpha}}{x^{0.5}}\lnc^{32}+\frac{p^\alpha}{x^{0.7}}\lnc^{32} +\frac{p^{1.5\alpha}}{x}\lnc^{33}\right)\right),
$$
where the constant in the $O$-notation depends on $\alpha$.}

We note that this formula becomes nontrivial if
$$
p^\alpha\ll x^{\frac{2}{3}}\lnc_x^{-\frac{68}{3}}.
$$

{\corollary\label{Sledstvie naimen chislo Hardy-Littlw} Let $q=p^\alpha$, where $p$ is a prime number and $(l,p)=1$. Then
$$
H_2(q,l)\ll q^{\frac{3}{2}}(\ln q)^{34}.
$$}

In proving the theorem \ref{Teorema o raspr chisel Hardy-Litlwoda}, we will use the results of T.~Cochrane \cite{Todd Cochrane-2002-Acta Arith-101} on estimating complete mixed sums $S(\chi,g,f,p^\beta)$, $\beta\ge2$. Note that the method for estimating complete sums of characters of the form $S(\chi,g,0,p^\beta)$ was developed by D.~Ismoilov~\cite{Ismoilov-DANRT-1986-10,Ismoilov-DANRT-1986-11,Ismoilov-DANRT-1990,Ismoilov-TrMIAN-1991,Ismoilov-ActaMathSin-1993,Ismoilov-TrMIAN-1994}, using the explicit formula of A.~G.~Postnikov~ \cite{PostnikovAG-1955-IzvANSSSR}.

{\sc Notation:}\vspace{-2pt}
\begin{itemize}
\item $x$~---~a sufficiently large positive real number;
\item $q$~~a natural number, $q>q_0$, $\chi(n)$~---~the Dirichlet character modulo $q$;
\item $\mu(n)$~---~the Mobius function; $s=\sigma +it$~--~a complex number; $M_j$, $N_j$, and $U_j$~--~integers, $N_j\le U_j<2N_j$;
\begin{align*}
 & S_j(s,\chi)=\sum\limits_{U_j< n\leq 2N_j}\dfrac{\chi(n)}{n^s}, \qquad G_j(s,\chi)=\sum\limits_{M_j< m\le 2M_j}\dfrac{\mu(m)\chi(m)}{m^s}; \\
 & W_k(s,\chi)=\sum\limits_{j=1}^kG_j(s,\chi)S_j(s,\chi), \qquad t_k(q;M,N)=\sideset{}{''}\sum\limits_{\chi} \int\limits_{-T}^T\left|W_k(0.5+it,\chi)\right|\dfrac{dt}{1+|t|};
\end{align*}
\item $\sideset{}{''}\sum\limits_\chi$~---~means summation over all primitive characters modulo $d$, $d|q$;
\item $ord_p(x)$~---~the greatest power of a prime number $p$ dividing an integer $x$, for a polynomial $f$ over $\mathbb{Z}$ $ord_p(f)$~---~the greatest power of $p$ dividing all coefficients of $f$, and for a rational function $f_1/f_2$, $ord_p(f_1/f_2)=ord_p(f_1)-ord_p(f_2)$.
\end{itemize}

\section{Known Lemmas}

{\lemma \label{Lemma teoremi o srednem} {\rm \cite{Montgomeri-1974}}. Assume that $M \geq 0$ and $N \geq 1$. Then, for any $T \geq 2$, the following inequality holds:
\begin{align*}
&\sideset{}{''}\sum_\chi\int^T_{-T}\left|\sum^{M+N}_{n=M+1}a_n\chi(n)n^{-it}\right|^2\frac{dt}{1+|t|}\ll \sum^{M+N}_{n=M}(n+q\ln T)|a_n|^2.
\end{align*}}

{\lemma \label{Lemma otsenka Mardzhan} {\rm \cite{Manjanishvili-1939-DANSSSr}}. For $x \geq 2$, we have:
$$
\sum_{n\leq x}\tau^2_{r}(n)\ll x(\ln x)^{r^2-1}.
$$}

{\lemma\label{Lemma todjestva Perrona} {\rm \cite{Karatsuba-OATCh}}. Let $b>0$ and $T>1$. Then, the following relation holds:
$$
\frac{1}{2\pi i}\int^{b+iT}_{b-iT}\frac{a^s}{s}ds=\left\{
\begin{array}{ll}
1+O\left(\dfrac{a^b}{T_0|\ln a|}\right), & \text{if} \quad \ a>1,  \\
O\left(\dfrac{a^b}{T_0|\ln a|}\right), & \text{if} \quad \ 0<a<1.
\end{array}
\right.
$$}

{\lemma \label{Lemma otsenka L-ryada v krit pol} {\rm \cite{Prakhar-1967-Mir}}. Let $q\geq 1$. Then, for $\operatorname{Re} s=\sigma\ge 0.5$, the following estimate holds:
$$
|L(s,\chi)|\ll (q|s|)^{1-\sigma}\ln q|s|.
$$}

{\lemma\label{Lemma momenti L-ryadov} {\rm \cite{Montgomeri-1974,Prakhar-1967-Mir}}. For $T\geq 2$, the following inequalities hold:
\begin{align*}
&\sideset{}{''}\sum_{\chi}\int^T_{-T}\left|L(0.5+it,\chi)\right|^4\frac{dt}{1+|t|}\ll q(\ln qT)^5,\\
&\sideset{}{''}\sum_{\chi}\int^T_{-T}\left|L(0.5+it,\chi)\right|^2\frac{dt}{1+|t|}\ll q\ln qT.
\end{align*}}

{\lemma\label{Lemma resheto} ( { \rm \cite{RakhmonovZKh-1993-IzvRAN}, using Heath-Brown identity~\cite{Heath-Brown-1982})}. Let $f(n)$ be an arbitrary complex-valued function, $u_1\leq x$, $r\geq 1$,
$$
C_r^k=\frac{r!}{k!(r-k)!},\quad \lambda (n)=\sum_{d\backslash n,\,d\le u_1}\mu(n).
$$
Then, the following identity holds:
\begin{align*}
\sum_{n\le x}\hspace{-2pt}\Lambda(n)f(n)=&\sum_{k=1}^r\hspace{-2pt}(-1)^{k-1}C_r^k\hspace{-4pt}\sum_{m_1\le u_1}\hspace{-6pt}\mu (m_1)\cdots\hspace{-20pt}\sum_{\substack{m_k\le u_1\\ m_1\cdots m_kn_1\cdots n_k\le x}}\hspace{-20pt}\mu (m_k) \sum_{n_1}\cdots \sum_{n_k}
                             \ln n_1f(m_1n_1\cdots m_km_k)\\
                            & +(-1)^r \sum_{n_1>u_1}\lambda (n_1)\cdots \sum_{\substack{n_r>u_1\\n_1\cdots n_rm\le x}}\lambda (n_r) \sum_{m}
                            \Lambda (m)f(n_1\cdots n_rm).
\end{align*}}

{\lemma\label{Lemma-otsen-poln-smshen-summ} {\rm \cite{Todd Cochrane-2002-Acta Arith-101}}. Let $f$ and $g$ be rational functions over $\mathbb{Z}$ that are not constant, let $p$ be an odd prime number, let $\chi$ be a Dirichlet character modulo $p^\beta$, and let $\beta$ be an integer such that $\beta \geq t+2$.

\textbf{(i)} If $\delta \not\in \mathcal{A}$, then $S_\delta(\chi,g,f,p^\beta)=0$.

If $\delta$ is a simple root, then
$$
S_\delta(\chi,g,f,p^\beta)=
\left\{
  \begin{array}{ll}
    \chi(g(\delta^*))e\left(\dfrac{f(\delta^*)}{p^\beta}\right)p^\frac{\beta+t}2, & \text{if } \beta - t \text{ is even,}  \vspace{4pt} \\
    \chi(g(\delta^*))e\left(\dfrac{f(\delta^*)}{p^\beta}\right)\left(\dfrac{A_\delta}p\right)\mathcal{G}_pp^\frac{\beta+t-1}2, & \text{if } \beta - t \text{ is odd,}
  \end{array}
\right.
$$
where $\delta^*$ is the comparison root constructed based on $\delta$,
\begin{align*}
\mathcal{C}(m) %:=p^{-t}(Rg(m)f'(x)+cg'(m))
\equiv0\left(\hspace{-9pt}\mod{p^{\left[\frac{\beta+t-1}2\right]}}\right),\qquad A_\delta\equiv2r(\mathcal{C}/g)'(\delta)\hspace{-9pt}\pmod p,
\end{align*}
where $\mathcal{G}$ is the quadratic Gauss sum.
% Given by the expression
%\begin{equation}\label{formula kv summa Gaussa}
%\mathcal{G}=\sum_{m=1}^pe_p(m^2)=\sum_{\delta=1}^{p-1}\left(\frac{m}{p}\right)e_p(m)=\left\{
%\begin{array}{ll}
%\sqrt{p}, & \text{if } p\equiv1\hspace{-9pt}\pmod p; \\
%i\sqrt{p}, & \text{if } p\equiv3\hspace{-9pt}\pmod p,
%\end{array}
% \right.
%\end{equation}
% and $\left(\frac{A_\delta}p\right)$ is the Legendre symbol.
}

\section{The Main Lemmas for Estimating Sums of the Form $t_k(q;M,N)$}

{\lemma\label{Lemma  otsenka tk(q;M,N) cherez teoremi o srednem} Let $T\geq 2$, $M_1\dots M_kN_1\dots N_k=Y$, $k_1\leq k$, $k_2\leq k$, $k_1+k_2=r$, $2^rM_{j_1}\dots M_{j_{k_1}}N_{i_1}\dots N_{i_{k_2}}=X$. Then, the following estimate holds:
\begin{align*}
&t_k(q;M,N)\ll (Y^{\frac12}+X^{\frac12}q^{\frac12}\lnc+Y^{\frac12}X^{-\frac12}q^{\frac12}\lnc+q\lnc^2)\lnc^{r^2+2k^2-2rk-1}.
\end{align*}}

{\sc Proof}. Let
\begin{align*}
&H_1(s,\chi)=\prod^{k_1}_{\alpha=1}G_{j_{\alpha}}(s,\chi)\prod^{k_2}_{t=1}S_{j_t}(s,\chi),\quad H_2(s,\chi)=W_k(s,\chi)H^{-1}_{1}(s,\chi).
\end{align*}
From the definitions of the functions $G_{j_{\alpha}}(s,\chi)$ and $S_{j_t}(s,\chi)$, it follows that
\begin{align*}
&H_1(s,\chi)=\sum_{n\le X}\frac{a_n\chi(n)}{n^s},\qquad H_2(s,\chi)=\sum_{m\le V}\frac{b_m\chi(m)}{m^s},\qquad V=2^{2k}YX^{-1},
\end{align*}
where $|a_{n}|\leq\tau_{r}(n)$, $|b_{m}|\leq\tau_{2k-r}(m)$. Applying Cauchy's inequality first for the integral over $t$ and then for the sum over characters $\chi$, we obtain:
\begin{align}
t_k(q;M,N)&=\sideset{}{''}\sum_{\chi}\int^T_{-T}\left|H_{1}(0.5+it,\chi)H_2(0.5+it,\chi)\right|\frac{dt}{1+|t|}\le \left(t_k'(q;M,N)t_k''(q;M,N)\right)^{\frac12},\label{formula otsenka tk(q;M,N)}\\
&t_k'(q;M,N)=\sideset{}{''}\sum_{\chi}\int^T_{-T}\left|\sum_{n\leq X}\frac{a_{n}\chi(n)}{n^{0,5+it}}\right|^2\frac{dt}{1+|t|},\nonumber\\
&t_k''(q;M,N)=\sideset{}{''}\sum_{\chi}\int^T_{-T}\left|\sum_{n\leq V}\frac{b_{n}\chi(n)}{n^{0,5+it}}\right|^2\frac{dt}{1+|t|}.\nonumber
\end{align}
Let's estimate $t_k'(q;M,N)$. Applying Lemmas \ref{Lemma teoremi o srednem} and \ref{Lemma otsenka Mardzhan} successively, we obtain
\begin{align*}
t_k'(q;M,N)&\ll \sum_{n\leq X}(n+q\ln T)\left|\frac{a_n}{n^{0,5}}\right|^2\ll\sum_{n\leq X}\tau^2_{r}(n)+q\ln T\sum_{n\leq X}\frac{\tau^2_{r}(n)}{n}\ll
%X\lnc^{r^2-1}+q\lnc^{r^2+1}\ll
(X+q\lnc^2)\lnc^{r^2-1}.
\end{align*}

In the same way we find that
\begin{align*}
t_k''(q;M,N)&\ll \sum_{n\leq V}(n+q\ln T)\left|\frac{b_n}{n^{0,5}}\right|^2\ll\left(\frac YX+q\lnc^2\right)\lnc^{(2k-r)^2-1}.
\end{align*}
From here and from the estimate of the sum  $t_k'(q;M,N)$, in view of (\ref{formula otsenka tk(q;M,N)}) the assertion of the lemma follows.

{\lemma \label{lemma virazhen S(0.5+it,chi) cherez L-ryad} Let $|t|\leq T$, $T\leq T_0$ and $N\le U<2N$, then the following inequality holds:
$$
S(0.5+it,\chi)\ll \int^{T_0}_{-T_0}|L(0.5+i(u+t),\chi)|\frac{du}{1+|u|}+\frac{N^{\frac12}\lnc}{T_0}+\left(\frac{q}{T_0}\right)^{\frac12}\lnc.
$$}

{\sc Proof.} Using Euler's formula, Lagrange's theorem on finite differences in the form $\sin \varphi=\sin \varphi-\sin 0=\varphi \cos\theta\varphi$, $0\le\theta\le 1$, as well as applying a trivial estimate, we obtain
\begin{align}
\left|\frac{(2N)^{iu}-U^{iu}}{u}\right|&=\frac{\sqrt{2-2\cos(u(\ln2N-\ln U))}}{|u|}=\frac{2|\sin(0.5u(\ln2N-\ln U))|}{|u|}=\nonumber\\
&=(\ln2N-\ln U)|\cos(0.5\theta u(\ln2N-\ln U))|\le\nonumber\\
&\le\min\left(\ln2N-\ln U,\ \frac{2}{|u|}\right) \le\min\left(\ln 2,\ \frac{2}{|u|}\right)\le\frac{2\ln 2+2}{1+|u|}. \label{formula otsenka ((2N)iu-Uiu)/u<<1/(1+|u|)}
\end{align}
Without loss of generality, we assume that $U$ and $2N$ are half-integer numbers. Applying Perron's identity (Lemma \ref{Lemma todjestva Perrona}) with $T=T_0$ and $b=0.5+(\ln2N)^{-1}$, we obtain
\begin{align}
S(0.5+it,\chi)&-\frac{1}{2\pi i}\int^{b+iT_0}_{b-iT_0}L(0.5+it+u,\chi)\frac{(2N)^u-U^u}udu\ll R_1(2N,T_0)+R_1(U,T_0),\label{formula S(0.5+it,chi)-J<<R_1(2N,T_0)+R_1(U,T_0)} \\
&R_1(N,T_0)=\sum^{\infty}_{n=1}\frac{1}{n^{0.5}}R\left(\frac{N}{n}\right)=\frac{1}{T_0}\sum^{\infty}_{n=1}\frac{1}{n^{0.5}}\left(\frac{N}{n}\right)^b \left|\ln \left(\frac{N}{n}\right)\right|^{-1}, \nonumber
\end{align}
where $N$ is one of the half-integer numbers $2N$ and $U$. The inequalities $\frac{N}{2}\geq n\ge 2N$ and $\left|\ln\left(\frac{N}{n}\right)\right|\ge \ln 2$ are equivalent. Therefore, considering the relation $n^{0.5+b}=n^{1+(\ln 2N)^{-1}}>n$, we obtain
\begin{align}
R_1(N,T_0)&=\frac{N^b}{T_0}\left(\sum_{\frac{N}{2}\geq n\ge 2N}\left(n^{0.5+b}\left|\ln\left(\frac{N}{n}\right)\right|\right)^{-1}+\sum_{\frac{N}{2}<n\le 2N-1}\left(n^{0.5+b}\left|\ln\left(\frac{N}{n}\right)\right|\right)^{-1}\right)\le \nonumber\\
&\le \frac{N^\frac12}{T_0}\left(\frac{1}{\ln 2}\sum_{\frac{N}{2}\geq n\ge 2N}\frac1{n^{1+(\ln 2N)^{-1}}}+\frac{2}{N}\sum_{\frac{N}{2}<n\le 2N-1}\left|\ln \left(\frac{N}{n}\right)\right|^{-1}\right).\label{formula otsenka R1(N,T0)}
\end{align}
Denoting the last two sums by $R_{11}$ and $R_{12}$, we estimate each separately. $R_{11}$ is a convergent numerical series, i.e., $R_{11}\ll 1$. In $R_{12}$, the summation variable $n$ takes integer values, starting from the integer $N_1$ to the integer $2N-1$, where
$$
N_1=\left\{
      \begin{array}{ll}
         \dfrac{N+0.5}{2}+1 , & \text{if $N-0.5$ is odd;} \\ \\
         \dfrac{N-0.5}{2}+1 , & \text{if $N-0.5$ is even.}
      \end{array}
    \right.
$$
By splitting the sum over $n$ into two parts and then using the equivalence of inequalities $N_1\le n\le N-0.5$ and $0\le N-0.5-n\le N-0.5-N_1$, we obtain
\begin{align*}
R_{12}&=-\sum_{n=N_1}^{N-0.5}\left(\ln\left(\frac{n}{N}\right)\right)^{-1}+\sum_{n=N+0.5}^{2N-1}\left(\ln\left(\frac{n}{N}\right)\right)^{-1}=\\
 &=\sum_{n=0}^{N-0.5-N_1}\left(-\ln\left(1-\frac{n+0.5}{N}\right)\right)^{-1}+\sum_{n=0}^{N-1.5}\left(\ln\left(1+\frac{n+0.5}{N}\right)\right)^{-1}.
\end{align*}

Next, for $0\le n\le N-0.5-N_1$ and $0\le n\le N-1.5$, using respectively the relations
\begin{align*}
&-\ln\left(1-\frac{n+0.5}{N}\right)=\sum_{k=1}^\infty\frac{1}{k}\left(\frac{n+0.5}{N}\right)^k>\frac{n+0.5}{N},\\
&\ln \left(1+\frac{n+0.5}{N}\right)=\sum_{k=1}^\infty\frac{(-1)^{k-1}}{k}\left(\frac{n+0.5}{N}\right)^k>\frac{n+0.5}{2N},
\end{align*}
we obtain
\begin{align*}
&R_{12}<\sum_{n=0}^{N-0.5-N_1}\frac{N}{n+1.5}+\sum_{n=0}^{N-1.5}\frac{2N}{n+0.5}\ll N\ln N.
\end{align*}
Substituting this estimate and the estimate for $R_{11}$ into formula (\ref{formula otsenka R1(N,T0)}), we get
$$
R_1(N,T_0)\ll \frac{N^{\frac12}}{T_0}\left(1+\frac{1}{N}\cdot N\ln N\right)\ll\frac{N^{\frac12}\lnc}{T_0}.
$$
Substituting this estimate into (\ref{formula S(0.5+it,chi)-J<<R_1(2N,T_0)+R_1(U,T_0)}), and then shifting the contour of integration in the integral $J$ to the line $\operatorname{Re}z=0$, we obtain
\begin{align*}
S(0.5+it,\chi)=&\frac{1}{2\pi i}\int^{T_0}_{-T_0}\frac{(2N)^{iu}-U^{iu}}{u}L(0.5+i(t+u),\chi)du+\\
&+\frac{1}{2\pi i}\int^{0}_{b}\frac{(2N)^{u-iT_0}-U^{u-iT_0}}{u-iT_0}L(0.5+u+i(t-T_0),\chi)du+\\
&+\frac{1}{2\pi i}\int^{b}_{0}\frac{(2N)^{u+iT_0}-U^{u+iT_0}}{u+iT_0}L(0.5+u+i(t+T_0),\chi)du+O\left(\frac{N^{\frac12}\lnc}{T_0}\right).
\end{align*}
Using inequality (\ref{formula otsenka ((2N)iu-Uiu)/u<<1/(1+|u|)}) for estimating the first integral and the estimate $|L(0.5+u+i(t+T_0),\chi)|\ll \left(qT_0\right)^{0.5-u}\ln qT_0$ (Lemma \ref{Lemma otsenka L-ryada v krit pol}) for the two remaining integrals, we obtain
\begin{align*}
|S(0.5+it,\chi)|&\ll\int^{T_0}_{-T_0}|L(0.5+i(t+u),\chi)|\frac{du}{1+|u|}+ \frac{\left(qT_0\right)^{0.5}\ln qT_0}{T_0}\int_0^b\left(\frac{N}{qT_0}\right)^udu+\frac{N^{\frac12}\lnc}{T_0}\ll\\
&\ll\int^{T_0}_{-T_0}|L(0.5+i(t+u),\chi)|\frac{du}{1+|u|}+\frac{N^{\frac12}\lnc}{T_0}+\left(\frac{q}{T_0}\right)^{\frac12}\lnc.
\end{align*}

{\lemma\label{lemma otsenka tk(q;M,N) cherez 4-moment} Let $ M_1\ldots  M_kN_1\ldots  N_k=Y$, $Y\leq y$, $y\leq x$, $N_1\geq N_2\geq\ldots  \geq N_k$, $ T\geq N_1^{\frac12}$, $q\leq T$. Then the following estimates hold:
\begin{align*}
&t_k(q;M,N)\ll\left(\left(\frac{Yq}{N_1N_2}\right)^{\frac12}+q\lnc\right)\lnc^{(2(k-1)^2+4)},\\
&t_k(q;M,N)\ll\left(\left(\frac{Yq}{N_1}\right)^{\frac12}+q\lnc\right)\lnc^{(2(k-0.5)^2+1)}.
\end{align*}}

{\sc Proof.} We have
\begin{align*}
&W_k(s,\chi)=H(s,\chi)S_1(s,\chi)S_2(s,\chi),\\
&H(s,\chi)=\prod_{j=1}^{k}G_j(s,\chi)\prod_{i=3}^kS_i(s,\chi)=\sum_{V_{0}<n\leq V_{1}}\frac{a(n)\chi(n)}{n^s}, \\
& a_n=\sum_{M_1<m_1\le 2M_1}\hspace{-10pt}\mu(m_1)\ldots \hspace{-10pt}\sum_{M_k<m_k\le 2M_k}\hspace{-13pt}\mu (m_k)\sum_{U_3<n_3\le 2N_3}\hspace{-8pt}\ldots \hspace{-8pt}\sum_{U_k<n_k\le 2N_k}\hspace{-8pt}1,\quad  |a_n|\le\tau_{2k-2}(n),\\
&V_0=M_1\ldots  M_kU_3\ldots  U_k,\quad V_1=2^{2k-2}M_1\ldots  M_kN_3\ldots  N_k\ll \frac{Y}{N_1N_2}.
\end{align*}
Applying Cauchy's inequality first to the integral over $t$, then to the sum over characters $\chi$, we obtain
\begin{align*}
&t_k(q;M,N)\le
\left(t_k'(q;M,N)t_k''(q;M,N)\right)^{\frac12}, \\
&t_k'(q;M,N)=\sideset{}{''}\sum_{\chi}\int_{-T}^T\left|\sum_{V_{0}<n\leq V_{1}}\frac{a(n)\chi(n)}{n^s}\right|^{2}\frac{dt}{1+|t|},\\
&t_k''(q;M,N)=\sideset{}{''}\sum_{\chi}\int_{-T}^T|S_1(0.5+it,\chi)S_2(0.5+it,\chi)|^2\frac{dt}{1+|t|}.
\end{align*}

Let us estimate $t_k'(q;M,N)$. Applying Lemma \ref{Lemma teoremi o srednem} sequentially, the relation $|a_n|\le\tau_{2k-2}(n)$, and then Lemma \ref{Lemma otsenka Mardzhan}, along with the relation $V_1\ll\dfrac{Y}{N_1N_2}$, we obtain
\begin{align*}
t_k'(q;M,N)&\ll\sum_{V_0<n\leq V_1}(n+q\ln T)\left|\frac{a_n}{n^{0.5}}\right|^2\ll\sum_{V_0<n\le V_1}\tau_{2k-2}^2(n)+q\ln T\sum_{V_0<n\leq V_1}\frac{\tau_{2k-2}^2(n)}{n}\ll \\
&\ll V_1(\ln V_1)^{(2k-2)^2-1}+q\ln T(\ln V_1)^{(2k-2)^2}\ll\left(\frac{Y}{N_1N_2}+q\lnc^2\right)\lnc^{(2k-2)^2-1}.
\end{align*}

Now, let us proceed to the estimation of $t_k''(q;M,N)$. Applying Lemma \ref{lemma virazhen S(0.5+it,chi) cherez L-ryad} to the sums \linebreak $S_1(0.5+it,\chi)$ and $S_2(0.5+it,\chi)$, with $T_0=T$, and noting that $T_0\geq \max(N_j^{\frac12},q)$, we obtain
$$
S_j(0.5+it,\chi)\ll\int^{T}_{-T}|L(0.5+i(u+t),\chi)|\frac{du}{1+|u|}+\lnc, \quad j=1,2.
$$
Utilizing this relation, the inequality $(a+b)^4\ll a^4+b^4$, and applying Cauchy's inequality twice to the inner integral over the variable $u$, followed by the symmetry of the repeated integral over the variables $u$ and $t$, we derive
\begin{align*}
t_k''(q;M,N)&\ll\sideset{}{''}\sum_{\chi}\int_{-T}^T\left|\int^{T}_{-T}|L(0.5+i(u+t),\chi)|\frac{du}{1+|u|}+\lnc\right|^4\frac{dt}{1+|t|}\ll\\ &\ll\sideset{}{''}\sum_{\chi}\int^{T}_{-T}\left(\int^{T}_{-T}|L(0.5+i(u+t),\chi)|\frac{du}{1+|u|}\right)^4\frac{dt}{1+|t|}+q\lnc^5,\\
&\ll\lnc^3\sideset{}{''}\sum_{\chi}\int^{T}_{-T}\int^{T}_{-T}|L(0.5+i(u+t),\chi)|^4\frac{du}{1+|u|}\frac{dt}{1+|t|}+q\lnc^5=\\
&=2\lnc^3\sideset{}{''}\sum_{\chi}\int_{|t|\le T}\int_{|t|\le |u|\le T}|L(0.5+i(u+t),\chi)|^4\frac{du}{1+|u|}\frac{dt}{1+|t|}+q\lnc^5.
\end{align*}
Since $|u|\geq|t|$, it follows that
$$
1+|u|\geq 1+\frac{|u|+|t|}{2}\geq 1+\frac{|u+t|}{2}\geq\frac12(1+|u+t|),
$$
therefore,
\begin{align*}
t_k''(q;M,N)&\ll\lnc^3\sideset{}{''}\sum_{\chi}\int_{|t|\leq T}\int_{|t|\leq |u|\leq T}\left|L(0.5+i(u+t),\chi)\right|^4\frac{du}{1+|u+t|}\ \frac{du}{1+|u|}+q\lnc^5=\\
&=\lnc^3\sideset{}{''}\sum_{\chi}\int_{|t|\leq T}\int_{|t|\leq |v-t|\leq T}\left|L(0.5+iv,\chi)\right|^4\frac{dv}{1+|v|}\ \frac{du}{1+|u|}+q\lnc^5=\ll\\
&\ll \lnc^4\sideset{}{''}\sum_{\chi}\int_{-2T}^{2T}\left|L(0.5+iv,\chi)\right|^4\frac{dv}{1+|v|}+q\lnc^5.
\end{align*}
Using Lemma \ref{Lemma momenti L-ryadov}, we obtain the assertion of the lemma. The second assertion of the lemma is proved similarly, but instead of the fourth moment of the Dirichlet $L$-series, its second moment is used.

\section{Proof of Theorem \ref{Teorema o srednem znach funk Chebisheva}}

Let $\chi_d$ be a primitive character modulo $d$, and $\chi$ an induced character from $\chi_d$ modulo $q$, where $d|q$. Then $\psi(y,\chi)=\psi(y,\chi_d)+O(\lnc^2),$ and hence,
\begin{align}\label{formula t(x,q)<<sum psi(y,chi)}
t(x;q)&=\sum_{\chi\neq\chi_0}\max_{y\leq x}|\psi(y,\chi)|+\psi(x,\chi_0)\ll\sideset{}{''}\sum_{\chi}\max_{y\leq x}|\psi(y,\chi)|+x+\varphi(q)\lnc^2.
\end{align}

Assuming in Lemma \ref{Lemma resheto}, $u=y^{\frac{1}{4}}$, $r=4$, and $f(n)=\chi(n)$, we obtain
\begin{equation}
\label{formula T(chi,nu)=sum_Tk(chi,nu)}
\psi(y,\chi)=\sum_{k=1}^4(-1)^kC_4^k\tilde{\psi}_k(y,\chi),
\end{equation}
$$
\tilde{\psi}_k(y,\chi)=\sum_{m_1\le u}\mu (m_1)\cdots \sum_{m_k\le u}\mu (m_k)\sum_{n_1}\cdots\hspace{-19pt}\sum_{m_1\cdots m_kn_1\cdots n_k\le y}\ln n_1\chi(m_1n_1\cdots m_kn_k).
$$
Dividing in $\tilde{\psi}_k(y,\chi)$ the ranges of each variable $m_1,\cdots ,m_k,n_1, \cdots ,n_k$ into at most $\lnc$ intervals of the form $M_j<m_j\le 2M_j$, $N_j<n_j\le 2N_j$, $j=1,2,\cdots ,k$, we obtain at most $\lnc^{2k}$ sums of the form
\begin{align*}
\hat{\psi}_k(y,\chi&,M,N)=\hspace{-10pt}\sum_{M_1<m_1\le 2M_1}\hspace{-10pt}\mu (m_1)\cdots\hspace{-10pt}\sum_{M_k<m_k\le 2M_k}\hspace{-10pt}\mu (m_k)\hspace{-10pt}\sum_{N_1<n_1\le 2N_1}\hspace{-10pt}\cdots\hspace{-10pt}\sum_{\substack{N_k<n_k\le 2N_k\\ m_1n_1\cdots m_kn_k\le y}}\hspace{-18pt}\chi (m_1n_1\cdots m_kn_k)\ln n_1=\\
=\int\limits_1^{2N_1}\hspace{-5pt}&\sum_{M_1<m_1\le 2M_1}\hspace{-10pt}\mu(m_1)\cdots\hspace{-15pt}\sum_{M_k<m_k\le 2M_k}\hspace{-13pt}\mu (m_k)\sum_{\max(u,N_1)<n_1\le 2N_1}\hspace{-5pt}\cdots\hspace{-5pt}\sum_{\substack{N_k<n_k\le 2N_k\\ m_1n_1\cdots m_kn_k\le x}}\hspace{-14pt}\chi (m_1n_1\cdots m_kn_k)d\ln u.
\end{align*}
Denoting by $U_1=\max (u,N_1)$ such a number $u$ at which the modulus of the integrand attains its maximum value, we obtain
\begin{equation}
\label{formula hatpsik(y,chi,M,N)<<lncpsik(y,chi,M,N)}
\hat{\psi}_k(y,\chi,M,N)|\ll \lnc \left|\psi_k(y,\chi,M,N)\right|,
\end{equation}
where
$$
\psi_k(y,\chi,M,N)=\sum_{M_1<m_1\le 2M_1}\hspace{-10pt}\mu(m_1)\cdots\hspace{-10pt}\sum_{M_k<m_k\le 2M_k}\hspace{-10pt}\mu (m_k)\sum_{U_1<n_1\le 2N_1}\cdots\sum_{\substack{U_k<n_k\le 2N_k\\ m_1n_1\cdots m_kn_k\le y}}\hspace{-10pt}\chi (m_1n_1\cdots m_kn_k),
$$
where $N_j\leq U_j<2N_j$, $j=1,2,\ldots ,k$. Without loss of generality, we assume that $M_1\ldots M_kN_1\ldots N_k<y$ and $y$ is a half-integer. The restriction $m_1n_1\ldots m_kn_k\le y$ is removed using Lemma \ref{Lemma todjestva Perrona} with $T=(xq)^{10}$:
\begin{align*}
\psi_k(y,&\chi,M,N)=\frac{1}{2\pi i}\int^{0.5+iT}_{0.5-iT}\prod_{j=1}^k\sum_{M_j<m_j\leq 2M_j}\frac{\chi(m_j)\mu(m_j)}{m_j^s}\sum_{U_j<n_j\leq 2N_j}\frac{\chi(n_j)}{n_j^s}\frac{y^s}{s}ds+\\
&\quad+O\left(\sum_{M_1<m_1\le 2M_1}\hspace{-10pt}m_1^{-\frac 12}\cdots\hspace{-10pt}\sum_{M_k<m_k\le 2M_k}\hspace{-10pt}m_k^{-\frac 12}\sum_{U_1<n_1\le 2N_1}\hspace{-10pt}n_1^{-\frac 12}\cdots\hspace{-10pt}\sum_{U_k<n_k\le 2N_k}\hspace{-8pt}n_k^{-\frac 12}
\frac{y^\frac{1}{2}}{T\left|\ln\frac{y}{m_1n_1\ldots m_kn_k}\right|}\right).
\end{align*}
For $m_1n_1...m_k n_k < y$ using the inequalities
\begin{align*}
\ln\frac{y}{m_1n_1\ldots m_kn_k}&\ge \ln\frac{y}{y-0.5}=\ln\left(1+\frac{1}{2y-1}\right)>\frac{1}{2y},
\end{align*}
and for $m_1n_1...m_k n_k > y$ using the inequalities
\begin{align*}
\ln\frac{m_1n_1\ldots m_kn_k}{y}&\ge \ln\frac{y+0.5}{y}=\ln\left(1+\frac{1}{2y}\right)>\frac{1}{2y},
\end{align*}
we obtain
\begin{align*}
\psi_k(y,&\chi,M,N)=\frac{1}{2\pi i}\int^{0.5+iT}_{0.5-iT}W_k(s,\chi)\frac{y^s}{s}ds+O\left(\frac{y^\frac{3}{2}}{T}\prod_{j=1}^k\sum_{M_j<m_j\leq 2M_j}m_j^{-\frac12}\sum_{N_j<n_j\leq 2N_j}n_j^{-\frac12}\right)=\\
&=\frac{1}{2\pi i}\int^{0.5+iT}_{0.5-iT}W_k(s,\chi)\frac{y^s}{s}ds+O\left(\frac{y^2}{(xq)^{10}}\right)\ll y^{\frac12}\int^T_{-T}|W_k(0.5+it,\chi)|\frac{dt}{1+|t|}+\frac{y^2}{(xq)^{10}}.
\end{align*}

Substituting the obtained estimate into (\ref{formula hatpsik(y,chi,M,N)<<lncpsik(y,chi,M,N)}), and then into (\ref{formula T(chi,nu)=sum_Tk(chi,nu)}), we obtain
\begin{align*}
\psi(y,\chi)%&\ll \sum_{k=1}^4|\tilde{\psi}_k(y,\chi)|\ll  \sum_{k=1}^4\lnc^{2k}|\tilde{\psi}_k(y,\chi,M,N)|\ll\\
&\ll y^{\frac12}\lnc^9\sum_{k=1}^4\int^T_{-T}|W_k(0.5+it,\chi)|\frac{dt}{1+|t|}+\frac{y^2\lnc^9}{(xq)^{10}}.
\end{align*}
From this and from the formula (\ref{formula t(x,q)<<sum psi(y,chi)}), we have
\begin{align}\label{formula t(x,q)<<sumtk(q;M,N)+..}
&t(x,q)\ll x^{\frac12}\lnc^9\sum_{k=1}^4\max_{y\le x}t_k(q;M,N)+x+\varphi(q)\lnc^2,\\
&t_k(q;M,N)=\sideset{}{''}\sum_{\chi_q}\int^T_{-T}|W_{k}(0.5+it,\chi)|\frac{dt}{1+|t|}. \nonumber
\end{align}
We estimate $t_k(q;M,N)$ separately for each $k=1,2,3,4$. Without loss of generality, we assume the following conditions for $t_k(q;M,N)$:
\begin{align}
&M_1\geq M_2\geq\ldots \geq M_k, \label{formula M1>M2>...MK}  \\
&N_1\geq N_2\geq\ldots \geq N_k, \label{formula N1>N2>...>NK} \\
&M_1\ldots M_kN_1\ldots  N_k=Y,   \quad Y \leq y,     \quad M_j\leq y^{\frac14}. \label{formula M1...MkN1...Nk=Y Y<y Mj<y^1/4}
\end{align}

\textbf{2. Estimate of $t_1(q;M,N)$.} Using the second statement of Lemma \ref{lemma otsenka tk(q;M,N) cherez 4-moment}, we get
\begin{align*}
&t_1(q;M,N)\ll\left(\left(M_1q\right)^{\frac12}+q\lnc\right)\lnc^{1,5}\ll(y^{\frac18}q^{\frac12}+q\lnc)\lnc^{1,5}\le (y^{\frac{3}{10}}q^{\frac12}+q\lnc)\lnc^{1,5}.
\end{align*}

\textbf{3. Estimate of $t_2(q;M,N)$.} Applying the first statement of Lemma \ref{lemma otsenka tk(q;M,N) cherez 4-moment}, we find
\begin{align*}
t_2(q;M,N)&\ll\left( \left( M_1M_2q\right)^{\frac12}+q\lnc\right)\lnc^6\le(y^{\frac{1}{4}}q^{\frac12}+q\lnc)\lnc^6\le (y^{\frac{3}{10}}q^{\frac12}+q\lnc)\lnc^6.
\end{align*}

\textbf{4. Estimate of $t_3(q;M,N)$.} Consider the three possible cases:
\begin{description}
\item[1]  $M_1M_2M_3\leq Y^{\frac25}$;
\item[2]  $Y^{\frac25}<M_1M_2M_3\leq Y^{\frac35}$;
\item[3]  $Y^{\frac35}<M_1M_2M_3$.
\end{description}

\textbf{Case 1. $M_1M_2M_3\leq Y^{\frac25}$}. Given the conditions (\ref{formula N1>N2>...>NK}) and (\ref{formula M1...MkN1...Nk=Y Y<y Mj<y^1/4}), we find that
\begin{align*}
&N_1N_2\ge N_1N_2\cdot\frac{N_3}{\sqrt[3]{N_1N_2N_3}}=\left(N_1N_2N_3\right)^\frac23=\left(\frac{Y}{M_1M_2M_3}\right)^\frac23\geq Y^{\frac25}.
\end{align*}
Therefore, by the first statement of Lemma \ref{lemma otsenka tk(q;M,N) cherez 4-moment}, we have
\begin{align*}
t_3(q;M,N)&\ll\left(\left(\frac{Yq}{N_1N_2}\right)^{\frac12}+q\lnc\right)\lnc^{12}\le \left(\left(Y^{\frac 35}q\right)^{\frac12}+q\lnc\right)\lnc^{12}\ll\left(y^{\frac{3}{10}}q^{\frac12}+q\lnc\right)\lnc^{12}.
\end{align*}

\textbf{Case 2. $Y^{\frac25}<M_1M_2M_3\leq Y^{\frac35}$.}  Using Lemma \ref{Lemma  otsenka tk(q;M,N) cherez teoremi o srednem} with $X=8M_1M_2M_3$, we get
\begin{align*}
t_3(q;M,N)&\ll\left(Y^\frac12+\left(M_1M_2M_3\right)^{\frac12}q^\frac12\lnc+Y^{\frac12}\left(M_1M_2M_3\right)^{-\frac12}q^{\frac12}\lnc+q\lnc^2\right)\lnc^8\le\\
&\leq\left(Y^{\frac12}+2Y^{\frac{3}{10}}q^{\frac12}\lnc+q\lnc^2\right)\lnc^8\ll\left(y^{\frac12}+y^{\frac{3}{10}}q^{\frac12}\lnc+q\lnc^2\right)\lnc^8.
\end{align*}

\textbf{Case 3. $Y^{\frac35}<M_1M_2M_3$.} From the relations (\ref{formula M1...MkN1...Nk=Y Y<y Mj<y^1/4}), (\ref{formula M1>M2>...MK}) and the condition of the case, we find that
\begin{align*}
&y^\frac12\ge M_1M_2\ge M_1M_2\cdot\frac{M_3}{\sqrt[3]{M_1M_2M_3}}=\left(M_1M_2M_3\right)^\frac{2}{3}\geq Y^{\frac25}.
\end{align*}
Therefore, by applying Lemma \ref{Lemma  otsenka tk(q;M,N) cherez teoremi o srednem} with $X=M_1M_2$, we have
\begin{align*}
t_3(q;M,N)&\ll\left(Y^{\frac12}+(M_1M_2)^{\frac12}q^{\frac12}\lnc+Y(M_1M_2)^{-\frac12}q^{\frac12}\lnc+q\lnc^2\right)\lnc^9\leq\\
&\le\left(Y^{\frac12}+Y^{\frac14}q^{\frac12}\lnc+Y^{\frac{3}{10}}q^{\frac12}\lnc+q\lnc^2\right)\lnc^9\ll
\left(y^{\frac12}+y^{\frac{3}{10}}q^{\frac12}\lnc+q\lnc^2\right)\lnc^9.
\end{align*}

\textbf{5. Estimate of $t_4(q;M,N)$.} Consider the seven possible cases:

\begin{description}
  \item[1] $M_1M_2M_3M_4\leq Y^{\frac15}$;
  \item[2] $Y^{\frac15}<M_1M_2M_3M_4\leq Y^{\frac25}$, $N_1N_2\leq Y^{\frac15}$;
  \item[3] $Y^{\frac25}<M_1M_2M_3M_4\leq Y^{\frac35}$, $N_1N_2\leq Y^{\frac15}$;
  \item[4] $Y^{\frac25}<M_1M_2M_3M_4\leq Y^{\frac35}$, $N_1N_2\geq Y^{\frac15}$;
  \item[5] $Y^{\frac35}<M_1M_2M_3M_4\leq Y^{\frac45}$, $N_1N_2\geq Y^{\frac15}$;
  \item[6] $Y^{\frac45}<M_1M_2M_3M_4\leq Y^{\frac55}$;
  \item[7] $Y^{\frac55}<M_1M_2M_3M_4$.
\end{description}

\textbf{Case 1. $M_1M_2M_3M_4 \leq Y^{\frac{1}{5}}$.} From the relations (\ref{formula N1>N2>...>NK}), (\ref{formula M1...MkN1...Nk=Y Y<y Mj<y^1/4}), and the conditions of the considered case, we have
$$
N_1N_2 \ge (N_1N_2N_3N_4)^\frac{1}{2} = \left(\frac{Y}{M_1M_2M_3M_4}\right)^\frac{1}{2} \ge Y^{\frac{2}{5}}.
$$
Therefore, according to the first statement of Lemma \ref{lemma otsenka tk(q;M,N) cherez 4-moment}, we obtain
\begin{align*}
t_4(q;M,N) & \ll \left( \left( \frac{Yq}{N_1N_2} \right)^\frac{1}{2} + q \lnc \right) \lnc^{22} \le \left( Y^{\frac{3}{10}} q^\frac{1}{2} + q \lnc \right) \lnc^{22} \leq \left( y^{\frac{3}{10}} q^\frac{1}{2} + q \lnc \right) \lnc^{22}.
\end{align*}

\textbf{Case 2. $Y^{\frac{1}{5}} < M_1M_2M_3M_4 \leq Y^{\frac{2}{5}}; \quad N_1N_2 \leq Y^{\frac{2}{5}}$.} From the relations (\ref{formula N1>N2>...>NK}), (\ref{formula M1...MkN1...Nk=Y Y<y Mj<y^1/4}), and the conditions of the considered case, we have
\begin{align*}
N_1N_2N_3 & \ge N_1N_2N_3 \frac{N_4}{\sqrt[4]{N_1N_2N_3N_4}} = (N_1N_2N_3N_4)^\frac{3}{4} \ge \left( \frac{Y}{M_1M_2M_3M_4} \right)^\frac{3}{4} \ge Y^\frac{9}{20} > Y^\frac{2}{5}, \\
N_1N_2N_3 & \le N_1N_2 \cdot \sqrt{N_1N_2} = (N_1N_2)^\frac{3}{2} \le \left( Y^\frac{2}{5} \right)^\frac{3}{2} = Y^\frac{3}{5}.
\end{align*}
Therefore, setting $X = 8N_1N_2N_3$ in Lemma \ref{Lemma otsenka tk(q;M,N) cherez teoremi o srednem}, we find
\begin{align*}
t_4(q;M,N) & \ll \left( Y^\frac{1}{2} + (N_1N_2N_3)^\frac{1}{2} q^\frac{1}{2} \lnc + Y^\frac{1}{2} (N_1N_2N_3)^{-\frac{1}{2}} q^\frac{1}{2} \lnc + q \lnc^2 \right) \lnc^{16} \leq \\
& \leq \left( Y^\frac{1}{2} + Y^\frac{3}{10} q^\frac{1}{2} \lnc + q \lnc^2 \right) \lnc^{16} \leq \left( y^\frac{1}{2} + y^\frac{3}{10} q^\frac{1}{2} \lnc + q \lnc^2 \right) \lnc^{16}.
\end{align*}

\textbf{Case 3. $Y^{\frac{1}{5}} < M_1M_2M_3M_4 \leq Y^{\frac{2}{5}}; \quad N_1N_2 > Y^{\frac{2}{5}}$.} Applying the first statement of Lemma \ref{lemma otsenka tk(q;M,N) cherez 4-moment} to the sum $t_4(q;M,N)$, we obtain
\begin{align*}
t_4(q;M,N) & \ll \left( \left( \frac{Yq}{N_1N_2} \right)^\frac{1}{2} + q \lnc \right) \lnc^{22} \le \left( Y^\frac{3}{10} q^\frac{1}{2} + q \lnc \right) \lnc^{22} \leq \left( y^\frac{3}{10} q^\frac{1}{2} + q \lnc \right) \lnc^{22}.
\end{align*}

\textbf{Case 4. $Y^{\frac{2}{5}} < M_1M_2M_3M_4 \leq Y^{\frac{3}{5}}$.} Applying Lemma \ref{Lemma otsenka tk(q;M,N) cherez teoremi o srednem} with $X = 16M_1M_2M_3M_4$, we have
\begin{align*}
t_4(q;N) & \ll \left( Y^\frac{1}{2} + Y^\frac{3}{10} q^\frac{1}{2} \lnc + q \lnc^2 \right) \lnc^{15} \leq \left( y^\frac{1}{2} + y^\frac{3}{10} q^\frac{1}{2} \lnc + q \lnc^2 \right) \lnc^{15}.
\end{align*}

\textbf{Case 5. $Y^{\frac{3}{5}} < M_1M_2M_3M_4 \leq Y^{\frac{4}{5}}; \quad M_1M_2M_3 \leq Y^{\frac{3}{5}}$.} From the relations (\ref{formula M1>M2>...MK}) and the conditions of the considered case, we find
\begin{align*}
M_1M_2M_3 & \ge M_1M_2M_3 \frac{M_4}{\sqrt[4]{M_1N_2M_3M_4}} = (M_1M_2M_3M_4)^\frac{3}{4} \ge \left( Y^\frac{3}{5} \right)^\frac{3}{4} = Y^\frac{9}{20} > Y^\frac{2}{5}.
\end{align*}
Therefore, using Lemma \ref{Lemma otsenka tk(q;M,N) cherez teoremi o srednem} with $X = 8M_1M_2M_3$, we find
\begin{align*}
t_4(q;M,N) & \ll \left( Y^\frac{1}{2} + Y^\frac{3}{10} q^\frac{1}{2} \lnc + q \lnc^2 \right) \lnc^{16} \le \left( y^\frac{1}{2} + y^\frac{3}{10} q^\frac{1}{2} \lnc + q \lnc^2 \right) \lnc^{16}.
\end{align*}

\textbf{Case 6. $Y^{\frac{3}{5}} < M_1M_2M_3M_4 \leq Y^{\frac{3}{5}}; \quad M_1M_2M_3 > Y^{\frac{3}{5}}$.} From the relations (\ref{formula M1...MkN1...Nk=Y Y<y Mj<y^1/4}), (\ref{formula M1>M2>...MK}), and the conditions of the considered case, we find
\begin{align*}
y^\frac{1}{2} & \ge M_1M_2 \ge M_1M_2 \frac{M_3}{\sqrt[3]{M_1M_2M_3}} = (M_1M_2M_3)^\frac{2}{3} > \left( Y^\frac{3}{5} \right)^\frac{2}{3} = Y^\frac{2}{5}.
\end{align*}
Therefore, in Lemma \ref{Lemma otsenka tk(q;M,N) cherez teoremi o srednem}, setting $X = 4M_1M_2$ and noting that $Y^\frac{2}{5} \ll X \ll y^\frac{1}{2}$, we have
\begin{align*}
t_4(q;M,N) & \ll \left( Y^\frac{1}{2} + y^\frac{1}{4} q^\frac{1}{2} \lnc + Y^\frac{3}{10} q^\frac{1}{2} \lnc + q \lnc^2 \right) \lnc^{19} \ll \left( y^\frac{1}{2} + y^\frac{3}{10} q^\frac{1}{2} \lnc + q \lnc^2 \right) \lnc^{19}.
\end{align*}

\textbf{Case 7. $M_1M_2M_3M_4 > Y^{\frac{4}{5}}$.} From the relations (\ref{formula M1...MkN1...Nk=Y Y<y Mj<y^1/4}), (\ref{formula M1>M2>...MK}), and the conditions of the considered case, we find
\begin{align*}
y^\frac{1}{2} & \ge M_1M_2 \ge M_1M_2 \frac{M_3M_4}{\sqrt{M_1M_2M_3M_4}} = (M_1M_2M_3M_4)^\frac{1}{2} > Y^\frac{2}{5}.
\end{align*}
Therefore, by Lemma \ref{Lemma otsenka tk(q;M,N) cherez teoremi o srednem} with $X = 4M_1M_2$, and noting that $Y^\frac{2}{5} \ll X \ll y^\frac{1}{2}$, we find:
\begin{align*}
t_4(q;M,N) & \ll \left( Y^\frac{1}{2} + Y^\frac{1}{4} q^\frac{1}{2} \lnc + Y^\frac{3}{10} q^\frac{1}{2} \lnc + q \lnc^2 \right) \lnc^{19} \ll \left( y^\frac{1}{2} + y^\frac{3}{10} q^\frac{1}{2} \lnc + q \lnc^2 \right) \lnc^{19}.
\end{align*}

\textbf{6.} Thus, for all $k = 1, 2, 3, 4$, it is proven that
\begin{align*}
\max_{y \leq x} t_k(q;M,N) & \ll \max_{y \leq x} \left( y^\frac{1}{2} \lnc^{19} + y^\frac{3}{10} q^\frac{1}{2} \lnc^{22} + q \lnc^{23} \right) = x^\frac{1}{2} \lnc^{19} + x^\frac{3}{10} q^\frac{1}{2} \lnc^{22} + q \lnc^{23}.
\end{align*}
Substituting these estimates into (\ref{formula t(x,q)<<sumtk(q;M,N)+..}), we obtain the statement of the theorem.

\section{Proof of Theorem \ref{Teorema otsenka S(a/q,x)} and its corollaries}
Using the orthogonality property of characters, we obtain
\begin{align*}
&S\left(\frac{a}q,x\right)=\frac1{\varphi(q)}\sum_{\chi\hspace{-7pt}\mod q}\chi(a)\tau(\bar{\chi})\psi(x,\chi)+O(\lnc^2),\\
&\tau(\chi)=\sum_{h=1}^q\chi(h)e\left(\frac hq\right),\qquad |\tau(\chi)|\ll\sqrt{q}.
\end{align*}
From this, using the relation $q\ll\varphi(q)\ln\lnc$ and Theorem \ref{Teorema o srednem znach funk Chebisheva}, we obtain
\begin{align*}
S\left(\frac{a}q,x\right)&\ll\frac{\sqrt{q}}{\varphi(q)}\sum_{\chi\hspace{-7pt}\mod q}\max_{y\le x}|\psi(y,\chi)|+\lnc^2 \ll\frac{\ln\lnc}{\sqrt{q}}t(x;q)+\lnc^2\ll \\
&\ll xq^{-\frac12}\lnc^{29}+x^{\frac{4}{5}}\lnc^{32}+x^\frac{1}{2}q^\frac12\lnc^{33}.
\end{align*}

To prove Corollary \ref{Sledstvie1 teoremi otsenka S(alpha,x)}, introducing the notation $\alpha-\frac aq=\lambda$, we consider two possible cases: $|\lambda|\le2x^{-1}$ and $2x^{-1}<|\lambda|\le q^{-2}$.

\textsc{Case 1: $|\lambda|\le2x^{-1}$.} Using Abel's transformation, we express the sum $S(\alpha,x)$ in terms of the sum $S\left(\frac aq, u\right)$, $u\le x$. We have
$$
S(\alpha,x)=-\int_{2}^{x}S\left(\frac aq, u\right)2\pi i\lambda e(u\lambda)du+e(\lambda x)S\left(\frac aq,x\right).
$$
Proceeding to the estimates and using the condition of the considered case, we find
\begin{align*}
&|S(\alpha,x)|\ll(|\lambda|x+1)\max_{u\le x}\left|S\left(\frac aq, u\right)\right|\ll xq^{-\frac12}\lnc^{29}+x^{\frac{4}{5}}\lnc^{32}+x^\frac{1}{2}q^\frac12\lnc^{33}.
\end{align*}

\textsc{Case 2: $2x^{-1}<|\lambda|\le q^{-2}$.} We have $\frac{q^2}{x}\leq \frac12$. According to Dirichlet's theorem on the approximation of real numbers by rational numbers, for any $\tau \ge 1$ there exist integers $b$ and $r$ that are coprime, with $1\le r\le \tau$, such that
\begin{align*}
&\left|\alpha-\frac br\right|\le\frac{1}{r \tau} .
\end{align*}
Taking $\tau=\frac xq$, we obtain
\begin{align}
&\left|\alpha-\frac br\right|\le\frac{q}{rx},\qquad r\le\frac xq. \label{formula br}
\end{align}	
Assume that $r=q$, then (\ref{formula br}) takes the form $|\lambda|\le\frac1x$, and as in Case 1, we obtain the required estimate for $S(\alpha,x)$. Now let $r\ne q$, then
\begin{align*}
&\left|\frac aq-\frac br\right|=\frac{|ar-bq|}{rq}\ge\frac{1}{rq}.
\end{align*}
From this and from $\frac{q^2}{x}\leq \frac12$, we obtain
\begin{align*}
\frac1{q^2}&\ge|\lambda|=\left|\left(\frac aq-\frac br\right)+\left(\frac br-\alpha\right)\right|\ge\left|\frac aq-\frac br\right|-\left|\alpha-\frac br\right|\ge\frac{1}{rq}-\frac{q}{rx}=\frac{1}{rq}\left(1-\frac{q^2}{x}\right)\ge\frac{1}{2rq},
\end{align*}
that is, $\frac qr\le2$. Therefore, (\ref{formula br}) takes the form:
\begin{align*}
&\left|\alpha-\frac br\right|\le\frac2x,\qquad \frac q2<r\le\frac xq.
\end{align*}
Consequently, as in Case 1, we obtain
\begin{align*}
&S(\alpha,x)\ll xr^{-\frac12}\lnc^{29}+x^{\frac45}\lnc^{32}+x^\frac12r^\frac12\lnc^{33}\ll xq^{-\frac12}\lnc^{33}+x^{\frac45}\lnc^{32}.
\end{align*}

Corollary \ref{Sledstvie2 teoremi otsenka S(alpha,x)} follows directly from Corollary \ref{Sledstvie1 teoremi otsenka S(alpha,x)}.

\section{Proof of Theorem \ref{Teorema o raspr chisel Hardy-Litlwoda}}
\quad \ Splitting the sum $\EuScript{H}_2(x;p^\alpha,l)$ into three parts and taking into account that $p^\alpha>\sqrt{x}$, we obtain
\begin{align*}
\EuScript{H}_2(x;p^\alpha,l)&=\sum_{\substack{n\leq x\\ (n,p)=1}}\Lambda(n)\sum_{\substack{m^2\leq x,\ (m^2-l,p)=1\\ n\equiv\,l-m^2\hspace{-7pt}\pmod{p^\alpha}}}1+\EuScript{H}_1'(x;p^\alpha,l)+\EuScript{H}_2''(x;p^\alpha,l),\\
\EuScript{H}_1'(x;p^\alpha,l)&=\sum_{\substack{n\leq x\\ (n,p)>1}}\Lambda(n)\sum_{\substack{m^2\le x \\ m^2\equiv\,l-n\hspace{-7pt}\pmod{p^\alpha}}}1\le  2\left(\frac{\sqrt{x}}{p^\alpha}+1\right)\lnc_x^2,\\
\EuScript{H}_2''(x;p^\alpha,l)&=\sum_{\substack{n\leq x\\ (n,p)=1}}\Lambda(n)\sum_{\substack{m\le\sqrt{x},\ m^2-l\equiv\,0\hspace{-7pt}\pmod p\\ m^2-l\equiv\,-n\hspace{-7pt}\pmod{p^\alpha}}}1=0.
\end{align*}
Next, using the orthogonality property of characters, we find
\begin{align*}
\EuScript{H}_2(x;p^\alpha,l)&=\frac{1}{\varphi(p^\alpha)}\sum_{\chi \hspace{-8pt}\mod{p^\alpha}}\psi(x,\chi) V_2(\sqrt{x},\overline{\chi},l,p^\alpha)+ O\left(\left(\frac{\sqrt{x}}{p^\alpha}+1\right)\lnc_x^2\right), \\
&V_2(u,\chi,l,p^\alpha)=\sum_{m\le u}\chi(l-m^2).
\end{align*}
Splitting the latter sum over $\chi$ into two parts, we obtain
\begin{align}
&\EuScript{H}_2(x;p^\alpha,l)=\EuScript G_2(x;p^\alpha,l)+\EuScript{R}_2(x;p^\alpha,l)+O\left(\left(\frac{\sqrt{x}}{p^\alpha}+1\right)\lnc_x^2\right), \label{formula H2=G2+R2+..}\\
&\EuScript{G}_2(x;p^\alpha,l)=\frac{\psi(x,\chi_0)V_2(\sqrt{x},\chi_0,l,p^\alpha)}{\varphi(p^\alpha)},\nonumber\\
&\EuScript{R}_2(x;p^\alpha,l)=\frac{1}{\varphi(p^\alpha)}\sum_{\chi\neq\chi_0}\psi(x,\chi)V_2(\sqrt{x},\overline{\chi},l,p^\alpha). \nonumber
\end{align}
In this formula, $\EuScript{G}_2(x;p^\alpha,l)$ gives the expected main term of $\EuScript{H}_2(x;p^\alpha,l)$, while $\EuScript{R}_2(x;p^\alpha,l)$ contributes to its remainder term.

We compute the main term. From the theorem of Ch.~Vallee-Poussin, we obtain
$$
\psi(x,\chi_0)=\sum_{n\leq x}\Lambda(n)+O(\lnc_x^{2})=x+O(x\exp(-c\sqrt{\lnc_x})).
$$
Now consider
\begin{align*}
V_2(\sqrt{x},\chi_0,l,p^\alpha)&=\sum_{m\le\sqrt{x}}1-\sum_{\substack{m\le\sqrt{x}\\(m^2-l,p)=p}}1 =\left[\sqrt{x}\right]-\sum_{\substack{m\le\sqrt{x}\\m^2\equiv l\hspace{-7pt}\pmod p}}\sum_{\substack{1\le n\le p\\n\equiv m\hspace{-7pt}\pmod p}}1=\\
&=\left[\sqrt{x}\right]-\sum_{\substack{1\le n\le p\\n^2\equiv l\hspace{-7pt}\pmod p}}\sum_{\substack{m\le\sqrt{x}\\m\equiv n\hspace{-7pt}\pmod p}}1 =\left[\sqrt{x}\right]-\sum_{\substack{1\le n\le p\\n^2\equiv l\hspace{-7pt}\pmod p}}\left[\frac{\sqrt{x}-n}p\right]=\\
&=\left[\sqrt{x}\right]-\sum_{\substack{1\le n\le p\\n^2\equiv l\hspace{-7pt}\pmod p}}\left(\frac{\sqrt{x}}p+O(1)\right) =x^\frac12\left(1-\frac{\rho(p,l)}{p}\right)+O(1),
\end{align*}
where $\rho(p,l)$ is the number of solutions of the congruence $n^2\equiv l\hspace{-5pt}\pmod p$, $1\le n\le p$. Therefore,
\begin{align}
\EuScript{G}_2(x;p^\alpha,l)&=\frac{x^{\frac32}}{\varphi(p^\alpha)}\left(1-\frac{\rho(p,l)}{p}+O\left(\exp(-c\sqrt{\lnc_x})\right)\right). \label{formula G2(x,palpha)}
\end{align}

We estimate the remainder term $\EuScript{R}_2(x;p^\alpha,l)$. Transitioning to primitive characters, we obtain
$$
\EuScript{R}_2(x;p^\alpha,l)=\frac{1}{\varphi(p^\alpha)}\sum_{\beta=1}^\alpha\ \sideset{}{^*} \sum_{\chi\hspace{-9pt}\mod\hspace{-2pt}p^\beta}\psi(x,\chi)V_2(\sqrt{x},\overline{\chi},l,p^\beta),
$$
where $^*$ indicates that summation is taken over primitive characters. Denoting by $\alpha_1$, where $1\le\alpha_1\le\alpha$, an integer satisfying the condition $p^{\alpha_1-1}\le\lnc_x^A<p^{\alpha_1}$, and then splitting the sum over $\beta$ into two parts $1\le\beta\le\alpha_1-1$ and $\alpha_1\le\beta\le\alpha$, we represent $\EuScript{R}_2(x;p^\alpha,l)$ as the sum of two terms $\EuScript{R}_{21}$ and $\EuScript{R}_{22}$. We first estimate $\EuScript{R}_{21}$. We have
\begin{align*}
\EuScript{R}_{21}&=\frac{1}{\varphi(p^\alpha)}\sum_{\beta=1}^{\alpha_1-1}\ \sideset{}{^*} \sum_{\chi\hspace{-9pt}\mod\hspace{-2pt}p^\beta}\psi(x,\chi)V_2(\sqrt{x},\overline{\chi},l,p^\beta)\ll\\
&\ll \frac{1}{\varphi(p^\alpha)}\sum_{\beta=1}^{\alpha_1-1}\ \sideset{}{^*}\max|\psi(x,\chi) | \ \sideset{}{^*}\sum_{\chi\hspace{-9pt}\mod\hspace{-2pt}p^\beta}|V_2(\sqrt{x},\chi,l,p^\beta)|,
\end{align*}
where the $^*$ in the sum over $\beta$ indicates that the maximum is taken over all primitive characters modulo $p^\beta$. Using the classical bound for $p^\beta\le\lnc_x^A$, $1\le\beta\le\alpha_1-1$ (see \cite{Karatsuba-OATCh}, p.~152),
$$
\psi(x,\chi)\ll x\exp\left(-c_1\sqrt{\lnc_x}\right),
$$
we obtain
\begin{align*}
\EuScript{R}_{21}&\ll\frac{x}{\varphi(p^\alpha)}\exp\left(-c_1\sqrt{\lnc_x}\right)\sum_{\beta=1}^{\alpha_1-1}\ \sideset{}{^*}\sum_{\chi\hspace{-7pt}\mod\hspace{-2pt}p^\beta}|V_2(\sqrt{x},\chi,l,p^\beta)|=\\
&=\frac{x}{\varphi(p^\alpha)}\exp\left(-c_1\sqrt{\lnc_x}\right)\sum_{\substack{\chi\hspace{-7pt}\mod\hspace{-2pt}p^{\alpha_1-1}\\ \chi\neq\chi_0}}|V_2(\sqrt{x},\chi,l,p^{\alpha_1-1})|.
\end{align*}

Next, applying Cauchy's inequality and then using the condition $p^{\alpha_1-1}\le\lnc_x^A$, we obtain
\begin{align}
\EuScript{R}_{21}&\ll\frac{x}{\varphi(p^\alpha)}\exp\left(-c_1\sqrt{\lnc_x}\right) \left(\varphi(p^{\alpha_1-1}) \sum_{\chi\hspace{-7pt}\mod\hspace{-2pt}p^{\alpha_1-1}}|V_2(\sqrt{x},\chi,l,p^\beta)|^2\right)^\frac12\ll\nonumber\\
&\ll \frac{x}{\varphi(p^\alpha)}\exp\left(-c_1\sqrt{\lnc_x}\right) \left(p^{\alpha_1-1}\sqrt{x}\left(\frac{\sqrt{x}}{p^{\alpha_1-1}}+1\right)\right)^\frac12\ll
\frac{x^{\frac32}}{\varphi(p^\alpha)}\exp\left(-c_1\sqrt{\lnc_x}\right).\label{formula R21}
\end{align}

Now we estimate $\EuScript{R}_{22}$. We have
\begin{align*}
\EuScript{R}_{22}&=\frac{1}{\varphi(p^\alpha)}\sum_{\beta=\alpha_1}^\alpha\ \sideset{}{^*} \sum_{\chi\hspace{-9pt}\mod\hspace{-2pt}p^\beta}\psi(x,\chi)V_2(\sqrt{x},\overline{\chi},l,p^\beta)\le \\
&\le\frac{1}{\varphi(p^\alpha)}\sum_{\beta=\alpha_1}^\alpha\ \sideset{}{^*}\max_{\chi\hspace{-9pt}\mod\hspace{-2pt}p^\beta} |V_2(\sqrt{x},\chi,l,p^\beta)|\sideset{}{^*}\sum_{\chi\hspace{-9pt}\mod\hspace{-2pt}p^\beta}|\psi(x,\chi)|,
\end{align*}
where the * symbol in the sum over $\beta$ indicates that the maximum is taken over all primitive characters modulo $p^\beta$. Using Theorem \ref{Teorema o srednem znach funk Chebisheva}, we obtain
\begin{equation}\label{formula R22(x,palpha,l)}
\EuScript{R}_{22}\ll\frac{x^\frac32}{\varphi(p^\alpha)}\sum_{\beta=\alpha_1}^\alpha\left(x^{-0.5}\lnc_x^{28}+x^{-0.7} p^\frac{\beta}2\lnc_x^{31}+ x^{-1}p^\beta\lnc_x^{32}\right)\ \sideset{}{^*}\max_{\chi\hspace{-9pt}\mod\hspace{-2pt}p^\beta}|V_2(\sqrt{x},\chi,l,p^\beta)|.
\end{equation}

We reduce the estimation of the incomplete sums $V_2(\sqrt{x},\chi,l,p^\beta)$ to the estimation of complete mixed sums of the form
$$
S(\chi,g,f,p^\beta)=\sum_{m=1}^{p^\beta}\chi(g(m))e\left(\frac{f(m)}{p^\beta}\right),\qquad g(m)=l-m^2, \quad f(m)=hm.
$$
We have the equality
\begin{align*}
&V_2(\sqrt{x},\chi,l,p^\beta)=\frac1{p^\beta}\sum_{h=1}^{p^\beta}\sum_{m\le \sqrt{x}}e\left(-\frac{hm}{p^\beta}\right)S(\chi,g,f,p^\beta)=\\
&=\frac{S(\chi,g,0,p^\beta)}{p^\beta}\left[\sqrt{x}\right]+\frac1{p^\beta}\sum_{h=1}^{p^\beta-1}\frac{\sin\frac{\pi h\sqrt{x}}{p^\beta}}{\sin\frac{\pi h}{p^\beta}}e\left(-\frac{h(1+\left[\sqrt{x}\right])}{2p^\beta}\right)S(\chi,g,f,p^\beta).
\end{align*}

Proceeding to estimates, we find
\begin{align*}
|V_2(\sqrt{x},\chi,l,p^\beta)|&\le\frac1{p^\beta}\max_{1\le h\le p^\beta}|S(\chi,g,f,p^\beta)|\left(\sqrt{x}+2\sum_{h=1}^{0.5(p^\beta-1)}\left(\sin\frac{\pi h}{p^\beta}\right)^{-1}\right).
\end{align*}
Since $p^\beta$ is an odd number, using the inequalities $\sin \pi \alpha \ge 2\alpha$ for $0\le\alpha<0.5$ and $\frac1h\le\ln\frac{2h+1}{2h-1}$ sequentially, we find
\begin{align*}
&2\sum_{h=1}^{0.5(p^\beta-1)}\left(\sin\frac{\pi h}{p^\beta}\right)^{-1} \le2\sum_{h=1}^{0.5(p^\beta-1)}\left(\frac{2h}{p^\beta}\right)^{-1} \le p^\beta\sum_{h=1}^{0.5(p^\beta-1)}\left(\ln(2h+1)-\ln(2h-1)\right)=p^\beta\ln p^\beta.
\end{align*}
Thus,
$$
|V_2(\sqrt{x},\chi,l,p^\beta)|\le\left(\frac{\sqrt{x}}{p^\beta}+\ln p^\beta\right)\max_{1\le h\le p^\beta}|S(\chi,g,f,p^\beta)|.
$$
Substituting this estimate into formula (\ref{formula R22(x,palpha,l)}), we obtain
\begin{equation}\label{formula otsenka R22}
\EuScript{R}_{22}\ll\frac{x^\frac32}{\varphi(p^\alpha)}\sum_{\beta=\alpha_1}^\alpha\left(\frac{\lnc_x^{28}}{x^{0.5}}+ \frac{p^\frac{\beta}2\lnc_x^{31}}{x^{0.7}}+\frac{p^\beta\lnc_x^{32}}{x}\right)\left(\frac{\sqrt{x}}{p^\beta}+\ln p^\beta\right)\sideset{}{^*}\max_{\substack{\chi\hspace{-9pt}\mod\hspace{-2pt}p^\beta,\\1\le h\le p^\beta}}\left|S(\chi,g,f,p^\beta)\right|,
\end{equation}
where the * symbol in the sum over $\beta$ indicates that the maximum is taken over all primitive characters modulo $p^\beta$. Next, we represent the sum $S(\chi,g,f,p^\beta)$ in the form
\begin{equation}\label{formulaS(chi,g,f,..)=sumdelta(chi,g,f,..)}
S(\chi,g,f,p^\beta)=\sum_{\delta=1}^pS_\delta,\qquad S_\delta=S_\delta(\chi,g,f,p^\beta)=\sum_{\substack{m=1\\ m\equiv\delta\hspace{-9pt}\pmod p}}^{p^\beta}\chi(g(m))e\left(\frac{f(m)}{p^\beta}\right).
\end{equation}

Let $a$ be the smallest primitive root modulo $p^\beta$. Define the number $r$ by the relation $a^{p-1}=1+rp$, $(r,p)=1$, and let $c=c(\chi,a)$ be the unique integer, $0<c\le p^{\beta-1}(p-1)$, such that for any integer $k$ the relation
$$
\chi(a^k)=e\left(\frac{ck}{ p^{\beta-1}(p-1)}\right)
$$
holds, i.e., the character $\chi$ is uniquely determined by the values of $r$ and $c$. Since in formula (\ref{formula otsenka R22}) all characters $\chi$ are primitive, it follows that $(c,p)=1$. Let
$$
t=t_p(\chi,g,f)=ord_p(rgf'+cg').
$$

Let $\mathcal{A}(\chi,g,f)$ denote the set of roots of the congruence
$$
\mathcal{C}(m):=p^{-t}(rg(m)f'(x)+cg'(m))\equiv0\hspace{-2pt}\pmod p,
$$
for which the terms in $S(\chi,g,f,p^\beta)$ are defined, that is,
$$
\mathcal{A}=\mathcal{A}(\chi,g,f):=\{\delta\in \mathbb{F}_p:\quad \mathcal{C}(\delta)\equiv0\hspace{-9pt}\pmod p,\quad g(\delta)\not\equiv0\hspace{-9pt}\pmod p\}.
$$
Now we define the set $\mathcal{A}$ in the case where $g=g(m)=l-m^2$ and $f=f(m)=hm$, depending on the parameters $h$, $c$, and $r$, taking into account that
$$
t_p(\chi,g,f)=ord_p(r(m^2-l)h+2cm)=\min\left(ord(rh),ord(2c),ord(lh)\right)=\min\left(ord(h),0\right)=0,
$$
we obtain that the set $\mathcal{A}$ has the form
\begin{equation}\label{formula C(alpha):=}
\mathcal{A}=\{\delta\in \mathbb{F}_p:\quad r(\delta^2-l)h+2c\delta\equiv0\hspace{-9pt}\pmod p,\quad \delta^2-l\not\equiv0\hspace{-9pt}\pmod p\},
\end{equation}
that is, $\mathcal{A}$ is the set of solutions to a quadratic congruence modulo $p$ and consists of at most two solutions. Let us consider two possible cases.

\textbf{1. Case $(h,p)=p$}. The quadratic congruence in (\ref{formula C(alpha):=}) reduces to a linear congruence of the form $2c\delta\equiv0\hspace{-2pt}\pmod p$, which has a single solution $\delta=p$.

\textbf{2. Case $(h,p)=1$}. Multiplying both sides of the congruence in (\ref{formula C(alpha):=}) by the number $rh$, where $(rh,p)=1$, and completing the square, we obtain
$$
\left(rh\delta+c\right)^2\equiv c^2+lr^2h^2\hspace{-9pt}\pmod p, \qquad 1\le\delta\le p-1.
$$
By the assumption of the theorem, $-l$ is a quadratic nonresidue, therefore, the right-hand side of the obtained quadratic congruence does not vanish modulo $p$, i.e.,
$$
c^2\not\equiv-lr^2h^2\hspace{-2pt}\pmod p.
$$
It follows that in (\ref{formula C(alpha):=}), the quadratic congruence
\begin{itemize}
  \item has no solutions if the number $c^2+lr^2h^2$ is a quadratic nonresidue;
  \item has two distinct solutions if the number $c^2+lr^2h^2$ is a quadratic residue.
\end{itemize}
Consequently, if the quadratic congruence in (\ref{formula C(alpha):=}) is solvable, then all roots are distinct, and there are at most two of them. Therefore, according to Lemma \ref{Lemma-otsen-poln-smshen-summ}, the right-hand side of (\ref{formulaS(chi,g,f,..)=sumdelta(chi,g,f,..)}) consists of at most two terms of the form $S_\delta(\chi,g,f,p^\beta)$, corresponding to these roots, for which the equality
$$
|S_\delta(\chi,g,f,p^\beta)|=p^{\frac{\beta}2}.
$$
holds. Substituting this estimate into formula (\ref{formula otsenka R22}), we obtain
\begin{align*}
\EuScript{R}_{22}&\ll\frac{x^\frac32}{\varphi(p^\alpha)}\left(\sum_{\beta=\alpha_1}^\alpha\left(\frac{\lnc_x^{28}}{p^{0.5\beta}}+ \frac{\lnc_x^{31}}{x^{0.2}}+\frac{p^{0.5\beta}}{x^{0.5}}\lnc_x^{32}+\frac{p^\beta}{x^{0.7}}\lnc_x^{32} +\frac{p^{1.5\beta}}{x}\lnc_x^{33}\right)\right)\ll\\
&\ll\frac{x^\frac32}{\varphi(p^\alpha)}\left(\frac{\lnc_x^{28}}{p^{0.5\alpha_1}}+\frac{\lnc_x^{32}}{x^{0.3}}+ \frac{p^{0.5\alpha}}{x^{0.5}}\lnc_x^{32}+\frac{p^\alpha}{x^{0.7}}\lnc_x^{32} +\frac{p^{1.5\alpha}}{x}\lnc_x^{33}\right).
\end{align*}
Further, using the choice of the number $\alpha_1$, we obtain
\begin{align*}
\EuScript{R}_{22}&\ll\frac{x^\frac32}{\varphi(p^\alpha)}\left(\lnc_x^{-0.5A+28}+ \frac{p^{0.5\alpha}}{x^{0.5}}\lnc^{32}+\frac{p^\alpha}{x^{0.7}}\lnc^{32} +\frac{p^{1.5\alpha}}{x}\lnc^{33}\right).
\end{align*}
From this, from (\ref{formula R21}) and (\ref{formula G2(x,palpha)}), in view of (\ref{formula H2=G2+R2+..}), the statement of Theorem \ref{Teorema o raspr chisel Hardy-Litlwoda} follows.

%%%%%%%%%%%%%%%%q%%%%%%%%%%%%%%%%%%%%%%%%%%%%%%%%%%%%%%%%%%%%%%%%%%%%%%%%%%%%%%%%%%%%%%%%%%%%%%%%%

\end{document}